\documentclass[12pt, reqno]{amsart}

\usepackage[utf8]{inputenc}
\usepackage[english]{babel} 

\usepackage{amssymb}
\usepackage{amsthm}  
\usepackage{amsmath} 
\usepackage{enumerate} 
\usepackage{fix-cm} 
\usepackage{enumitem}
\usepackage[mathscr]{eucal} 
\usepackage{comment}

\calclayout

\usepackage[inner = 25mm, outer = 25mm, top = 25mm, bottom = 32mm, a4paper]{geometry}

\usepackage{tikz-cd}

\usepackage{aliascnt}

\usepackage{graphicx}

\makeatletter
\def\l@subsection{\@tocline{2}{0pt}{2.5pc}{5pc}{}}
\makeatother

\DeclareRobustCommand{\SkipTocEntry}[5]{}

\theoremstyle{definition}
\newtheorem{dfn}{Definition}[subsection]

\theoremstyle{remark}

\newaliascnt{rmk}{dfn}
\newtheorem{rmk}[rmk]{Remark}
\aliascntresetthe{rmk}

\newaliascnt{ex}{dfn}

\aliascntresetthe{ex}

\theoremstyle{plain}
\newaliascnt{thm}{dfn}
\newtheorem{thm}[thm]{Theorem}
\aliascntresetthe{thm}

\newaliascnt{prop}{dfn}
\newtheorem{prop}[prop]{Proposition}
\aliascntresetthe{prop}

\newaliascnt{prop2}{dfn}

\aliascntresetthe{prop2}

\newaliascnt{lem}{dfn}
\newtheorem{lem}[lem]{Lemma}
\aliascntresetthe{lem}

\newaliascnt{cor}{dfn}

\aliascntresetthe{cor}

\usepackage{xcolor} 
\usepackage[colorlinks=true, linkcolor=blue, urlcolor=red, citecolor=brown]{hyperref} 

\addto\extrasenglish{

  }

\let\oldfootnotemark\footnotemark
\let\oldfootnotetext\footnotetext
\let\oldfootnote\footnote
\renewcommand\footnote[1]{\addtocounter{footnote}{1}\hypertarget{fnbackref.\arabic{footnote}}{}\addtocounter{footnote}{-1}\oldfootnote{#1\fnbackref}}
\renewcommand\footnotemark{\addtocounter{footnote}{1}\hypertarget{fnbackref.\arabic{footnote}}{}\addtocounter{footnote}{-1}\oldfootnotemark}
\renewcommand\footnotetext[1]{\oldfootnotetext{#1\fnbackref}}
\newcommand{\fnbackref}{\hyperlink{fnbackref.\arabic{footnote}}{\footnotesize$\uparrow$}}

\DeclareMathOperator{\Hom}{Hom}

\DeclareMathOperator{\id}{id}
\DeclareMathOperator{\Tot}{Tot}

\newcommand{\bounded}{\mathrm{b}}

\newcommand{\mb}[1]{\mathbb{#1}}
\newcommand{\ca}[1]{\mathcal{#1}}
\newcommand{\mcr}[1]{\mathscr{#1}}

\newcommand\into{\hookrightarrow}
\newcommand{\pgap}{\vspace{5pt}}

\newcommand{\mD}{\mathrm{D}}

\newcommand{\mH}{\mathrm{H}}

\newcommand{\mR}{\mathrm{R}}

\newcommand{\cA}{\mathcal{A}}

\newcommand{\cD}{\mathcal{D}}
\newcommand{\cE}{\mathcal{E}}
\newcommand{\cF}{\mathcal{F}}
\newcommand{\cG}{\mathcal{G}}
\newcommand{\cH}{\mathcal{H}}
\newcommand{\cI}{\mathcal{I}}

\newcommand{\cO}{\mathcal{O}}

\newcommand{\C}{\mathbb{C}}

\newcommand{\bS}{\mathbb{S}}

\renewcommand\P{\mathbb{P}}
\newcommand\A{\mathbb{A}}
\newcommand\Z{\mathbb{Z}}
\newcommand\bbS{\mathbb{S}}

\newcommand{\rR}{\mathscr{R}}

\newcommand\Crit{\operatorname{Crit}}

\DeclareMathOperator{\MF}{MF}
\DeclareMathOperator{\ind}{ind}

\makeatletter
\newcommand{\nodeequation}[1]{%
  \let\label\ltx@label
  \refstepcounter{equation}%
  #1
  \quad
  (\theequation)%
}
\makeatother

\begin{document}


\title{Serre functors of residual categories via hybrid models}


\author{Federico Barbacovi}
\address{Department of Mathematics, University College London}
\email{federico.barbacovi.18@ucl.ac.uk}

\author{Ed Segal}
\address{Department of Mathematics, University College London}
\email{e.segal@ucl.ac.uk}

\begin{abstract} In this short note we observe that the Serre functor on the residual category of a complete intersection can be easily described in the framework of hybrid models. Using this description we recover some recent results of Kuznetsov and Perry.
\end{abstract}

\maketitle

\tableofcontents

\section{Introduction}\label{sec.intro}

Let $X\subset \P^n$ be a smooth Fano complete intersection of multi-degree $(d_1,..., d_k)$, and write $d=\sum_{i} d_i$. The ambient $\P^n$ has a very simple derived category since we have Beilinson's full exceptional collection consisting of the line bundles $\cO, ..., \cO(n)$. Kuznetsov \cite{KuzHPD} observed that we can use this fact to simplify the derived category of $X$ itself: if we take the subset $\cO,..., \cO(n-d)$ and restrict them to $X$ then they still form an exceptional collection, but no longer a full one. So they have some orthogonal $\rR_X$ and we get a semi-orthogonal decomposition:
\begin{equation}\label{eqn.SOD}\mD^{\bounded}(X) = \big\langle \rR_X, \, \cO, ..., \cO(n-d) \big\rangle \end{equation}
This is only a non-trivial decomposition because $X$ is Fano, \emph{i.e.}$d\leq n$; in the Calabi-Yau or general type case we have $\rR_X=\mD^{\bounded}(X)$. 

The subcategory $\rR_X$ is called the \emph{residual component} (or \emph{Kuznetsov component}) of $\mD^{\bounded}(X)$. It controls much of the interesting behaviour of $\mD^{\bounded}(X)$: the deformation theory, moduli spaces of objects, stability conditions, \emph{etc.} \cite{PS, BLMS, FV}.  It also has a Serre functor $\bbS$,  and one can ask about the behaviour of $\bbS$. For example in the hypersurface case Kuznetsov showed that $\rR_X$ is fractional Calabi-Yau, meaning that some power of $\bbS$ is a shift functor \cite{KuzCY}. 

More recently Kuznetsov and Perry \cite{KP} studied $\bbS$ for higher-codimension $X$ by a recursive procedure. They think of $X$ as the end point of a sequence 
\begin{equation}\label{eqn:slices}\P^n\supset X_{d_1} \supset X_{d_1, d_2}\supset ... \supset X_{d_1,.., d_k}=X\end{equation}
given by repeatedly intersecting with hypersurfaces, and relate the Serre functors on residual categories at each step. Using this result they were able to compute some quantities called the upper and lower \emph{Serre dimensions} of $\rR_X$ (see Section \ref{sec.Serredim}). This result corrects a conjecture of Kontsevich and Katzarkov. They also have some very nice applications, such as a proof of the non-existence of a Serre invariant stability condition on $\rR_X$ in most cases. 
\pgap

The purpose of this note is to give an alternative (and rather shorter) derivation of Kuznetsov and Perry's results, using an idea from the string theory literature called a `hybrid model'.

In the case when $X$ is a hypersurface there is a celebrated theorem of Orlov \cite{OrlovCohSing} which can be interpreted as relating $\rR_X$ to a category of matrix factorizations on an affine orbifold:
\begin{equation}\label{eq.Orlov}\rR_X = \MF\!\big( [\A^{n+1} / \Z_{d}], \; f \big) \end{equation}
Here $f$ is the degree $d$ polynomial defining $X$, re-interpreted as a $\Z_d$-invariant function on $\A^{n+1}$. This construction was generalized to complete intersections in \cite{Seg} (following \cite{witten}) using the following two steps:
\begin{enumerate}
\item 
Let $Y_+$ be the total space of the vector bundle $\bigoplus_{i=0}^k \cO(-d_i) $ over $\P^n$. On this we have a function:
$$ W= \sum_{i=1}^k f_i p_i$$
Here $f_1,.., f_k$ are the defining polynomials of $X$, pulled-up to $Y_+$, and each $p_i$ is the tautological section of $\cO(-d_i)$ on $Y_+$. Using Kn\"orrer periodicity (\emph{e.g.} \cite{OrlovKP, Shipman}) we equate $\mD^{\bounded}(X)$ with the category of matrix factorizations:
$$\mD^{\bounded}(X) \cong \MF\left( Y_+, W\right) $$

\item Let $\P(\mathbf{d})$ denote the weighted projective space $\P(d_1,..., d_k)$, and let $Y_-$ be the total space of the vector bundle $\cO(-1)^{n+1}$ over  $\P(\mathbf{d})$. The spaces $Y_-$ and $Y_+$ are related by an `orbifold flip', and there is a fully-faithful functor:
\begin{equation}\label{eqn.ff}\MF (Y_-, W ) \into \MF (Y_+,  W)\end{equation}
Note that on $Y_-$ the $p_i$ are the homogeneous co-ordinates of weighted projective space, and the $f_i$ are polynomials in the fibre co-ordinates. Under Kn\"orrer periodicity the image of this embedding is identified with $\rR_X$.
\end{enumerate}

The pair $(Y_-, W)$ can be viewed as a family, indexed over $\P(\mathbf{d})$, of the kind of affine-orbifolds-plus-a-function that appeared in the hypersurface case \eqref{eq.Orlov}. The latter are called \emph{Landau-Ginzburg models} in the physics literature,\footnote{Although this term is now used in the mathematics literature with a more general meaning.} and this `family of Landau-Ginzburg models' is called a \emph{hybrid model}.
\pgap

For us, the advantage of identifying $\rR_X$ with $\MF(Y_-, W)$ is that it makes the Serre functor very transparent; it is simply given (up to a shift) by tensoring with the canonical bundle $\omega_{Y_-}$. This follows from a general fact about Serre duality for matrix factorizations, which was proven in the generality we need by Favero and Kelly \cite{FavKel}. We recall it as Proposition \ref{prop:SD-MF} below. 
\pgap

 When $X$ is a hypersurface the fact that $\rR_X$ is fractional Calabi-Yau follows easily from this fact and Orlov's theorem \eqref{eq.Orlov}. Indeed on the orbifold $[\A^n/\Z_d]$ the $d$th power of any line bundle is trivial, so some power of $\bbS$ is a shift. 
  
When $X$ has higher codimension the properties of $\bbS$ are more complicated, but our point-of-view makes them fairly easy to prove. 
 
 \subsection*{Plan of the paper}
 \begin{itemize}\setlength{\itemsep}{10pt}
     
 \item In Section \ref{sec.background} we set up some general terminology for matrix factorizations, in particular we recall the notion of \emph{R-charge} which is essential for getting the correct shifts. 

We then produce an explicit generator for the category $\MF(Y_-, W)$, by pushing forward the usual exceptional collection of line bundles from $\P(\mathbf{d})$. The reason this works is that the critical locus of $W$ is some thickening of $\P(\mathbf{d})$, and it's enough to generate $\mD^{\bounded}(\Crit(W))$ since the category of matrix factorizations is supported there.
 
 \item  In Section \ref{sec.Serredim} we discuss Serre dimensions. Roughly, these are given by taking a generator for your category, applying the Serre functor repeatedly, and seeing how the maximal and minimal degrees of Ext groups grow. 
 
 
Since we have an explicit generator and a simple description of $\bbS$ we can find the Serre dimensions by a straight-forward calculation. 

\item In Section \ref{sec.twist} we explain how one more result of Kuznetsov-Perry can be understood in our framework. As explained above, they take a recursive approach where they repeatedly slice by hyperplanes and see how the Serre functor changes. In fact they identify some power of the Serre functor (up to a shift) with some power of a spherical twist induced by the previous slice \cite[Cor.~1.4]{KP}. 

 In our approach to Serre dimension we don't need this result, but it's still interesting to translate it into the hybrid model description. It turns out that their spherical twist is just a line bundle on $Y_-$, and then the observation is that some power of this bundle equals some power of $\omega_{Y_-}$.
 
 \end{itemize}
 
 \subsection*{Acknowledgements}
 
This project has received funding from the European Research Council (ERC) under the European Union Horizon 2020 research and innovation programme (grant agreement No.725010).

We thank the referee for pointing us to the reference \cite{FavKel}.

\section{Background}\label{sec.background}

\subsection{General theory}

Let $Y$ be a smooth variety or stack equipped with:
\begin{enumerate}
\item[(i)] A $\C^*$ action such that $-1\in \C^*$ acts trivially.
\item[(ii)] A regular function $W\in \Gamma(\cO_Y)$ of weight 2.
\end{enumerate}
We call (i) the \emph{R-charge} and (ii) the \emph{superpotential}. Given this data, a \emph{matrix factorization} is an equivariant sheaf $\cE$ on $Y$ equipped with a map $d: \cE\to \cE$ of R-charge 1, satisfying $d^2 = W \cdot \id_\cE$. Matrix factorizations form the objects of a  dg-category $\MF(Y,W)$ (\cite{Seg, OrlovNonaffine, BFK} \emph{etc.}). 

There is an equivariant line bundle on $Y$ associated to the generating character of $\C^*$, we denote it by $\cO[1]$. Tensoring by this line bundle is the shift operator on $\MF(Y,W)$. Note that in this formulation the category is $\Z$-graded, not $\Z_2$-periodic or $\Z_2$-graded.

\begin{rmk} \label{localExts}

In the category of matrix factorizations there is a derived local hom's functor which sends matrix factorizations  $\cE, \cF$ to an object $\mR\ca{H}om(\cE, \cF)\in \MF(Y, 0)$. To compute it we must either replace $\cE$ by a `locally free resolution', \emph{i.e.}~an equivalent matrix factorization whose underlying sheaf is a vector bundle, or $\cF$ by an `injective resolution'.   The derived global sections of $\mR\ca{H}om(\cE, \cF)$ give the morphisms from $\cE$ to $\cF$.

 More generally, if $\cE\in \MF(Y, W)$ but $\cF\in \MF(Y, W')$ for some other superpotential $W'$ then $\mR\ca{H}om(\cE, \cF)$ is an object in $\MF(Y, W'-W).$
\end{rmk}

\begin{prop} \cite[Theorem~2.18]{FavKel}
    \label{prop:SD-MF} Assume that the critical locus $\Crit(W)\subset Y$ is proper. Then the category $\MF(Y,W)$ is smooth and proper, and it has a Serre functor given by $- \otimes \omega_Y [\dim Y]$ where $\omega_Y$ is the equivariant canonical bundle.
\end{prop}

In fact Favero-Kelly work on a more general class of Artin quotient stacks. For earlier related results see \cite[Sect.~8]{preygel}, \cite{LP}, \cite[Lem.~6.8]{dyckerhoff}, \emph{etc.}

\begin{lem}
    \label{lem:generator-MF}
    Let $i : \Crit(W)_{rd} \hookrightarrow Y$ denote the reduced scheme underlying the critical locus of $W$. Let $\cG$ be a $\C^*$-equivariant sheaf on $\Crit(W)_{rd}$ which, up to twisting by all characters of $\C^*$, generates the equivariant derived category $\mD^{\bounded}_{\mathbb{C}^*}(\Crit(W)_{rd})$.    Then $i_{\ast}\cG$ is a generator of  $\MF(Y,W)$.
\end{lem}
\begin{proof}
    This follows from, for example, Proposition 3.64 and Corollary 4.14 of \cite{BFK}.
\end{proof}

\begin{rmk}The condition that $\cG$ is a sheaf (rather than a complex) is not essential but it means we can ignore the subtle distinction between $\mD^{\bounded}_{\mathbb{C}^*}(\Crit(W)_{rd})$ and $MF(\Crit(W)_{rd}, 0)$.\end{rmk}

\subsection{Our setup}

Now we specialize the general theory above to the example we're interested in. 
\pgap

 We fix positive integers $n, d_1,...,d_k$ and let $d= \sum_{i=1}^k d_i$. Let $\C^*$ act on $\C^{n+k+1}$ with weights $(d_1,..., d_k,
 -1,..., -1)$, and write $p_1,...,p_k, x_0,..., x_n$ for the corresponding co-ordinates.   We define an orbifolds $Y_\pm$ by taking the stack-theoretic quotients of two different open subsets:
  $$ Y_+ = \big [ \{(x_0,..., x_n)\neq 0\}\;/\; \C^*\big] $$
  $$ Y_- = \big [ \{(p_1,..., p_k)\neq 0\}\;/\; \C^*\big] $$
  $Y_+$ is the total space of a vector bundle over $\P^n$, and $Y_-$ is the total space of a vector bundle over the weighted projective space $\P(\mathbf{d})$.
  
  Now we introduce an R-charge on both these spaces by adding a second $\C^*$ action, which acts with weight 2 on each $p_i$ and weight 0 on each $x_i$. On $Y_+$ this action just rotates the fibres. On $Y_-$ the action  preserves the zero section $\P(\mathbf{d})\subset Y_-$ but does not act trivially on it, unless all the $d_i$'s are equal. 
  
  Finally we define the superpotential by choosing polynomials $f_1,...,f_k\in \mathbb{C}[x_0,...,x_n]$ with $f_i$ of degree $d_i$, then setting:
  $$ W= \sum_{i=1}^k p_i f_i $$
 With this data we have categories of matrix factorizations $\MF(Y_+, W)$ and $\MF(Y_-, W)$. As discussed in the introduction, Kn\"orrer periodicity \cite{Shipman} gives an equivalence
$$MF(Y_+, W)\cong \mD^{\bounded}(X)$$
 where $X\subset \P^n$ is the complete intersection cut out by the $f_i$. Moreover, if $d<n$ then we have an embedding
 $\MF(Y_-, W) \into \MF(Y_+, W)$
 whose image is the residual category $\rR_X$ (this follows from \cite{Seg}).
 
 We remark that neither of the two results above require that $X$ is smooth. However, we will make this assumption from now on. Given this, the following is an elementary calculation:
 \begin{lem} \label{lem:crit}The critical locus of $W$ on $Y_+$ is exactly $X$. On $Y_-$, the critical locus is a thickening of the zero section $\P(\mathbf{d})$.
  \end{lem}
  
 It now follows from Proposition \ref{prop:SD-MF} that $\MF(Y_-,W)$ has a Serre functor given by a shift of the canonical bundle. There's some potential confusion here due to R-charge, so we spell this out explicitly.
 
 Any equivariant line bundle on $Y_-$ is specified by its weights under the two $\C^*$ actions, so is of the form $\cO(k)[m]$. Recall here that we use square brackets for the R-charge.\footnote{The R-charge \emph{is} the shift in the category of matrix factorizations, there isn't a separate homological shift.} In this notation the equivariant canonical bundle of $Y_-$ is
 $$\omega_{Y_-} = \cO(n+1-d)[-2k] $$
 and hence the Serre functor is:
 \begin{equation}\label{eq:Serre}
     \omega_{Y_-}[\dim Y_-] = \cO(n+1-d)[n-k]
 \end{equation}

On $Y_+$ we get the same answer apart from a sign change, so there the Serre functor is $\cO(d-n-1)[n-k]$. Note that this agrees with the Serre functor on $X$, as required by Kn\"orrer periodicity. 
 
 We can also use Lemmas \ref{lem:generator-MF} and \ref{lem:crit} to deduce: 
 
 \begin{lem}\label{lem:split-generator}
Let $T$ be the (R-charge equivariant) vector bundle
$$ T = \bigoplus_{i=0}^{d-1} \cO(i) $$
on $\mb{P}(\mathbf{d})$. Then the push-forward of $T$ generates $\MF(Y_-,W)$.
\end{lem}

\section{Serre dimensions}\label{sec.Serredim}

Let $\mcr{D}$ be a triangulated category with a Serre functor $\bbS$. For objects $E,F\in \mcr{D}$ we write $\Hom_{\cD}^\bullet(E, F)$ for the graded space of morphisms between $E$ and all shifts of $F$. For simplicity assume that the Homs are bounded in degree.

The Serre dimensions of $\cD$, defined by Elagin and Lunts, are some kind of measure of the complexity of $\bbS$.

\begin{dfn}\cite{EL}\; 
\begin{enumerate}[leftmargin=*, itemsep=10pt]
    \item  For two objects $E,F \in \mcr{D}$ we define
  \[
    \begin{array}{lcr}
      e_{-}(E,F) = \min \{k : \Hom_{\mcr{D}}^k(E,F) \neq 0 \} & \mathrm{and} & e_{+}(E,F) = \max \{k : \Hom_{\mcr{D}}^k(E,F) \neq 0 \}.
    \end{array}
  \]
  
  \item
  Choose a generator $G\in \cD$. The \textbf{upper} and \textbf{lower Serre dimension} of $\cD$ are  $$ \overline{\mathrm{Sdim}}(\cD) = \limsup_{m\to \infty} \frac{-e_-(G, \bbS^m G)}{m} \quad\mbox{and}\quad \underline{\mathrm{Sdim}}(\cD) = \liminf_{m\to \infty} \frac{-e_+(G, \bbS^m G)}{m}.$$
  \end{enumerate}
\end{dfn}

As the notation suggests, these Serre dimensions do not depend on the choice of generator $G$.

\vspace{5pt}

We want to compute the Serre dimensions of the category $\rR_X\cong \MF(Y_-, W)$. We have a simple description of the Serre functor \eqref{eq:Serre} and we have an explicit generator (Lemma \ref{lem:split-generator}) so now we just need to carry out the computation.

For simplicity we assume that each $d_i\geq 2$, of course this can always be achieved by reducing the number of variables.

\begin{lem}\label{lem:epm} Let $d_{\max}$ and $d_{\min}$ denote the maximum and minimum of the $d_i$'s.  Then for $l\gg  0$ we have
  $$      e_{-}\big(\ca{O}_{\mb{P}(\mathbf{d})}, \ca{O}_{\mb{P}(\mathbf{d})}(l)\big) \;\geq\; \frac{2l}{d_{\max}} $$
and the bound is attained whenever $l$ is a multiple of $d_{\max}$. We also have
  $$   e_{+}\big(\ca{O}_{\mb{P}(\mathbf{d})}, \ca{O}_{\mb{P}(\mathbf{d})}(l)\big)\;\leq \;\frac{2(l-n-1)}{d_{\min}} + n+ 1 $$
 and the bound is attained whenever $l-n-1$ is a multiple of $d_{\min}$.   

\end{lem}

\begin{proof}
    We need to compute the morphisms between $\cO_{\P(\mathbf{d})}$ and $\cO_{\P(\mathbf{d})}(l)$, and we do this by first computing the local morphisms via a resolution of $\cO_{\P(\mathbf{d})}$ (see Remark \ref{localExts}).
  
  To get the resolution we observe that $\mb{P}(\mathbf{d}) \subset X$ is the zero locus of a regular section of $\cO_{Y_-}(-1)^{\oplus n + 1}$. We can write our superpotential in the form $W=\sum_{i} x_i g_i$ for some $g_i$'s, and then $\cO_{\P(\mathbf{d})}$ is equivalent to the `Koszul type' matrix factorization built from the $x_i$'s and the $g_i$'s; see \emph{e.g.}~\cite[Sect.~3.2]{Seg}. Our assumption on degrees ensures that each $g_i$ vanishes along $\P(\mathbf{d})$ so it follows that 
    \begin{equation}
    \label{eqn:ext-alg-rhom-tot-space}
    \mR \ca{H}om\big(\ca{O}_{\mb{P}(\mathbf{d})}, \ca{O}_{\mb{P}(\mathbf{d})}(l)\big)  =  \bigoplus_{i=0}^{n+1} \ca{O}_{\mb{P}(\mathbf{d})}(l-i)^{\oplus \binom{n+1}{i}}[-i]
  \end{equation}
  with no differential.  If $l \gg 0$ then these line bundles have only zeroth cohomology, and the morphisms from  $\cO_{\P(\mathbf{d})}$ to $\cO_{\P(\mathbf{d})}(l)$ are simply the global sections of the RHS of \eqref{eqn:ext-alg-rhom-tot-space}.
\pgap  

  We know $\mH^0(\ca{O}_{\mb{P}(\mathbf{d})}(j))$ is the degree $j$ part of $\mb{C}[p_1, \dots, p_k]$.  The question we must answer is: what are the minimum and maximum possible R-charges of a polynomial with degree $j$?   Notice that if $p_1^{a_1} \dots p_k^{a_k}$ has degree $j$, then $\sum a_i d_i = j$, and its R-charge is $2\sum a_i$. But then:
    \[
    2\frac{j}{d_{\max}} = 2\sum a_i \frac{d_i}{d_{\max}} \leq 2 \sum a_i
  \]
  Hence the minimum R-charge attainable is at least $2j/d_{\max}$.
  Taking into account the shift in \eqref{eqn:ext-alg-rhom-tot-space}, we see that 
  $$e_{-}\big(\ca{O}_{\mb{P}(\mathbf{d})}, \ca{O}_{\mb{P}(\mathbf{d})}(l)\big) \geq 
  \min_{i\in [0,n+1]}    \left( 2 \frac{(l-i)}{d_{\max}} + i \right) = 2 \frac{l}{d_{\max}}
  $$
  where we used $d_{\max} > 1$. Moreover this minimum degree is actually attained whenever $l$ is a multiple of $d_{\max}$, by a polynomial $p_t^{l/d_{\max}}$ where $p_t$ has degree $d_{\max}$. 
  
Similarly
  $$e_{+}\big(\ca{O}_{\mb{P}(\mathbf{d})}, \ca{O}_{\mb{P}(\mathbf{d})}(l)\big) \leq 
  \max_{i\in [0,n+1]}    \left( 2 \,\frac{(l-i)}{d_{\min}} + i \right) = 2 \frac{(l-n-1)}{d_{\min}} + n+1 $$
 and the bound is attained whenever $d_{\min}$ divides $l-n-1$. 
\end{proof}

Recall that the index of the Fano $X$ is $\ind X = n + 1 - d$. 

\begin{thm} \cite[Thm.~1.7]{KP}
  We have
  \[
    \begin{array}{lcr}
      \displaystyle{\overline{\mathrm{Sdim}}(\mcr{R}_X)= \dim X - 2\, \frac{\ind X}{d_{\max}}}
      & \mathrm{and} &
      \displaystyle{\underline{\mathrm{Sdim}}(\mcr{R}_X) = \dim X - 2\, \frac{\ind X}{d_{\min}}}
    \end{array}
  \]
\end{thm}

\begin{proof}
  We prove the first statement, the second being identical.
  By \eqref{eq:Serre} the Serre functor for $\MF(Y_-,W)$ is given by $- \otimes \ca{O}_{Y_-}(\ind X)[\dim X]$.
  Using the generator of \autoref{lem:split-generator}, we have to compute:
  \[
    \begin{aligned} 
       & \limsup_{m\to \infty}\; \frac{-1}{m}\,
        e_{-}\!\left(\bigoplus_{i=0}^{d-1} \ca{O}_{\mb{P}(\mathbf{d})}\big(i\big),\; \bigoplus_{j=0}^{d-1} \ca{O}_{\mb{P}(\mathbf{d})}\big(j + m\ind X\big)[m\dim X]\right) \\
= & \dim X  - \liminf_{m\to \infty}\,  \min_{j\in [1-d, d-1]}\, \frac{1}{m}\, e_{-}\!\Big( \ca{O}_{\mb{P}(\mathbf{d})},\;  \ca{O}_{\mb{P}(\mathbf{d})}\big(j + m\ind X\big)\Big)\\   
    \end{aligned}
\]
But Lemma \ref{lem:epm} implies that
$$ \frac{2(d-1 + m \ind X)}{m d_{\max}} \; \geq \min_{j\in [1-d, d-1]} \, \frac{1}{m}\, e_{-}\!\Big( \ca{O}_{\mb{P}(\mathbf{d})},\;  \ca{O}_{\mb{P}(\mathbf{d})}\big(j + m\ind X\big)\Big) \; \geq\;  \frac{2(1-d + m \ind X)}{m d_{\max}} $$ 
where for the upper bound, we observe that there some $j$ in this range such that $d_{\max}$ divides $j+m\ind X$. The result follows.

      
\end{proof}

\section{The recursive approach}\label{sec.twist}

As discussed in the introduction, Kuznetsov and Perry take a recursive approach, viewing $X$ as the end point of a sequence of complete intersections of increasing codimension. In this section we explain what hybrid models can contribute to this point-of-view.

 Let $M\subset \P^n$ be the penultimate variety in this sequence, \emph{i.e.} the variety denoted by $X_{d_1,..., d_{k-1}}$ in \eqref{eqn:slices}. So $X$ is a divisor in $M$, obtained by intersecting $M$ with a hypersurface of degree $d_k$. By construction the restriction functor $\mD^{\bounded}(M) \to \mD^{\bounded}(X)$ induces a functor
\begin{equation}\label{eq:restrict}\rR_M \to \rR_X.\end{equation}
Kuznetsov and Perry observe that this functor is spherical, and that the twist and cotwist around it are closely related to the Serre functors on $\rR_M$ and $\rR_X$. More precisely a certain power of the twist (resp.~cotwist) is equal, up to a shift, to a certain power of the Serre functor on $\rR_X$ (resp.~$\rR_M$). They deduce this as an instance of a more general result involving spherical functors and semi-orthogonal decompositions.

\vspace{5pt}
As before we view $\rR_X$ as a category of matrix factorizations on the orbifold $Y_-$. In the same way $\rR_M$ is a category of matrix factorizations on the orbifold
$$Z_{-} = \Tot \big\{\cO(-1)^{n+1} \to \P(d_1,..., d_{k-1})\big\}. $$
Notice the reversal of inclusions here: $X$ is a divisor in $M$, but $Z_-$ is a divisor in $Y_-$. Moreover the appropriate superpotential on $Z_-$ is just the restriction of $W$ from $Y_-$. This implies, writing
$$j: Z_- \into Y_-$$
for the inclusion, that we have a push-forward functor 
\begin{equation}\label{eq:pushforward}j_*: \MF(Z_-, W) \to \MF(Y_-, W).\end{equation}
It is not hard to verify that this agrees with the functor \eqref{eq:restrict}, but we won't do it here since we have not discussed the details of the equivalence $\MF(Y_-, W)\cong \rR_X$.

For ordinary derived categories it is a well-known and easy-to-prove fact that the push-forward along the inclusion of a divisor $j: Z\into Y$ induces a spherical functor $j_*: \mD^{\bounded}(Z)\to \mD^{\bounded}(Y)$ (e.g.~\cite[Sect.~1.2]{addington}). The twist around $j_*$ is given by tensoring with $\cO_Y(Z)$, and the cotwist by tensoring with $\cO_{Z}(Z)[-2]$.

The proofs of these facts extend to categories of matrix factorizations immediately, the only modification being that we replace homological shifts by R-charge equivariance. Since $Z_-\subset Y_-$ is the set $\{p_{d_k}=0\}$, and $p_{d_k}$ is a section of the equivariant line bundle $\cO_{Y_-}(d_k)[2]$, we can conclude immediately that the functor \eqref{eq:pushforward} is spherical with twist and cotwist:
\begin{equation}\label{eq:twistandcotwist} T_{j_*} = \otimes\, \cO_{Y_-}(d_k)[2], \hspace{3cm} C_{j_*} = \otimes\, \cO_{Z_-}(d_k) \end{equation}

Since the Serre functors on $\MF(Y_-, W)$ and $\MF(Z_-, W)$ are also line bundles, and both these varieties have Picard rank 1, we have proven:

\begin{prop}\cite[Cor.~1.4]{KP} Some power of $T_{j_*}$ agrees, up to a shift, with some power of the Serre functor on $\MF(Y_-, W)$. Similarly some power of $C_{j_*}$ agrees, up to a shift, with some power of the Serre functor on $\MF(Z_-, W)$.
\end{prop}
In fact comparing \eqref{eq:twistandcotwist} with \eqref{eq:Serre} we can compute that

$$\mb{S}_{\MF(Z_-,W)}^{d_k / c} \simeq C_{j_{\ast}}^{\ind(M)/c} \left[ \frac{d_k \dim(M)}{c} \right]\;\;\;\mbox{and}\;\;\;
						\mb{S}_{\MF(Y_-,W)}^{d_k / c} \simeq T_{j_{\ast}}^{\ind(X)/c} \left[ \frac{d_k \dim(X) - 2 \ind(X)}{c} \right]
		$$
where $c = \gcd(d_k, \ind(M))=\gcd(d_k, \ind(X))$. This is exactly Kuznetsov and Perry's result.

\bibliographystyle{halphanum}

\end{document}